\documentclass{article}
\usepackage[utf8]{inputenc}

\usepackage{amsmath}
\usepackage{amssymb}
\usepackage{amsfonts}
\usepackage{amsthm}
\usepackage{a4wide}
\usepackage{hyperref}
\usepackage{color}
\usepackage[english]{babel}
\usepackage{graphicx}
\usepackage[all]{xy}
\usepackage{blindtext}
\usepackage{xcolor}
\usepackage{authblk}

\theoremstyle{theorem}
\newtheorem{theorem}{Theorem}
\newtheorem{proposition}[theorem]{Proposition}
\newtheorem{lemma}[theorem]{Lemma}
\newtheorem{corollary}[theorem]{Corollary}

\theoremstyle{definition}

\newcommand{\D}{\mathcal{D}}

\newcommand{\Z}{\mathbb{Z}}
\newcommand{\R}{\mathbb{R}}
\newcommand{\Q}{\mathbb{Q}}
\newcommand{\K}{\mathbb{K}}

\begin{document}

\title{A number system with base $-\frac32$}
\markright{A number system with base $-\frac32$}
\author{Luc\'ia Rossi\footnote{The doctoral position of this author is supported by the Austrian Science Fund (FWF), project W1230.}$\;$ and J\"org M. Thuswaldner}
\affil{Montanuniversit\"at Leoben, Franz Josef Stra\ss{}e~18, A-8700 Leoben, Austria}
\date{}



\maketitle

\begin{abstract}
In the present paper we explore a way to represent numbers with respect to the base $-\frac32$ using the set of digits $\{0,1,2\}$. Although this number system shares several properties with the classical decimal system, it shows remarkable differences and reveals interesting new features. For instance, it is related to the field of $2$-adic numbers, and to some ``fractal'' set that gives rise to a tiling of a non-Euclidean space.
\end{abstract}

\section{Introduction.}

A number system is, intuitively, a way of representing a certain set of numbers in a consistent manner, using strings of some given digits in relation to a base. The most famous examples are the decimal and the binary systems. Over time, many generalizations of these number systems came to the fore. They have applications in  various areas of mathematics and computer science. To cite some examples, back in 1885 Gr\"unwald \cite{Gru} studied number systems with negative integers as bases. 
In 1936 
Kemp\-ner~\cite{Kempner} and later also R\'enyi \cite{Renyi:57} proposed expansions w.r.t.\ nonintegral real bases (see {\it e.g.} also \cite{MR2024754}). Knuth~\cite{Knuth:60} introduced complex bases and dealt with their relations to fractal sets (see \cite[p. 608]{Knu} for an important example). Gordon~\cite{MR1613189} discussed the relevance of number systems with varying digit sets in cryptography and in \cite{MR2042759} binary and hexadecimal expansions were used in this context.

\begin{figure}[h!]
			\centering\includegraphics[scale=0.6]{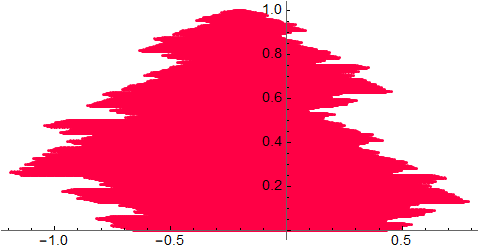}
			\caption{The tile $\mathcal{F}$ related to the number system with base $-\frac32$.}
			\label{tile}
		\end{figure}

The aim of this article is to define and explore a number system which has base $-\frac32$ and uses the set of digits $\{0,1,2\}$. We chose this example because it has many beautiful properties that link different areas of  mathematics, and nevertheless can be studied in a way that is accessible to a broad readership. The results stated here are known in a more general context; however, since our proofs are presented in terms of our particular example, we can avoid the use of advanced techniques. Number systems of a very similar kind were introduced by Akiyama {\it et al.}~\cite{MR2448050} and generalizations of this have been studied in \cite{MR2735753} and \cite{SSTW}. We mention that there are other ways to define a number system with base $-\frac32$, see for example the one using the digit set $\{0,1\}$ in~\cite{MR3051451}.

Proceeding in analogy to the decimal system, we first investigate expansions of integers and real numbers. Realizing that the real line is somehow ``too small'' for our number system, we introduce the space $\K=\R\times\Q_2$, where $\Q_2$ is the field of $2$-adic numbers, and show how it naturally arises as a representation space. We relate the set $\mathcal{F}$ depicted in Figure~\ref{tile} to our number system. This set has a self-affine structure and many nice properties. Among other things, we prove that $\mathcal{F}$ induces a tiling of $\K$ by translations, and relate this tiling to the existence and uniqueness of expansions w.r.t. the base $-\frac32$ in $\K$. The set $\mathcal{F}$ is an example of a \textit{rational self-affine tile}. Steiner and Thuswaldner~\cite{MR3391902} dealt with such tiles in a much more general framework.


\section{Expansions of the integers --- and more.} \label{section 2}

In the decimal system, each integer can be expanded without using digits after the decimal point. In this section, we wish to define and explore such ``integer expansions'' in the number system with base $-\frac32$ and digit set $\mathcal{D}=\{0,1,2\}$, which we denote as $(-\frac32,\D)$. 
The fact that $-\frac32$ is not an integer will entail new properties. The reason for the choice of a negative base is that there will be no need for a minus sign to represent negative numbers. 
		
To get a feeling for this number system, we first deal with the integers and follow the ideas of  \cite{MR2448050}. Given $N\in\mathbb{Z}$, we want to find an expansion of the form	
\begin{equation}\label{eq:Nrep}
N=\tfrac{1}{2}\sum_{i=0}^{k}d_i\left(-\tfrac32\right)^i \qquad (k\in\mathbb{N},\; d_i\in\mathcal{D}).
\end{equation}
 From here onwards, we assume $d_k\neq0$ whenever $k\geqslant 1$. The factor $\frac12$ at the beginning is just there for convenience (and to be consistent with~\cite{MR2448050}) and will not be crucial. In order to produce an expansion of the form (\ref{eq:Nrep}), we use the following algorithm. Write 
\begin{equation*}
2N=-3N_1+d_0
\end{equation*}
with $d_0\in\mathcal{D}$, $N_1\in\Z$. Since $N$ is given, $d_0$ has to be the unique digit satisfying $d_0\equiv 2N \pmod{3}$, so this equation has a unique solution $N_1$.
More generally, we set $N_0=N$ and recursively define the integers $N_{i+1}$ for $i\geqslant0$ by the equations 
\begin{equation}\label{eq:algo}
2N_i=-3N_{i+1}+d_i
\end{equation}
with $d_i\in \mathcal{D}$. One can easily see by induction on $i$ that this yields 
\begin{equation}\label{eq:NNi}
N=\left( -\tfrac{3	}{2}\right)^{i+1}N_{i+1}+\tfrac{d_{i}}{2}\left( -\tfrac{3}{2}\right)^{i}+\cdots+\tfrac{d_1}{2}\left( -\tfrac{3}{2}\right)+ \tfrac{d_0}{2}.
\end{equation}
If we can prove that $N_k=0$ for $k$ large enough, our algorithm gives the desired representation \eqref{eq:Nrep} for each $N\in \Z$. And this is indeed our first result.	 
	 
\begin{proposition}\label{prop:N}
Each $N\in \Z$ can be represented in the form \eqref{eq:Nrep} in a unique way.
\end{proposition}
	
\begin{proof}
By \eqref{eq:NNi}, it suffices to show that for each $N=N_0\in\mathbb{Z}$ the sequence $(N_i)_{i\geqslant0}$ produced by the recurrence \eqref{eq:algo} is eventually zero. We have $N_{i+1}=-\tfrac23N_i+\tfrac13d_i$ with $d_i\in\D$, hence
$|N_{i+1}|\leqslant\tfrac23|N_i|+\tfrac23$		
and therefore $|N_{i+1}| < |N_i|$ holds for each $|N_i| \geqslant 3$. This implies that there is $i\in\mathbb{N}$ with $|N_i| \leqslant 2$. 

Direct calculation shows that if $N_i=-2$ then $N_{i+1}=2$, $N_{i+2}=-1$, $N_{i+3}=1$, $N_{i+4}=0$, and $N_{i+5}=0$. 
Thus for each $N=N_0\in\mathbb{Z}$ there is $k_0\in \mathbb{N}$ with $N_{k}=0$ for all $k\geqslant k_0$ and $N$ has an expansion of the form \eqref{eq:Nrep}. 

Concerning uniqueness of the expansion, we just note that each digit $d_i$ in \eqref{eq:Nrep} has to lie in a prescribed residue class modulo $3$, and hence is uniquely determined.
\end{proof}

We write  	
\[
(d_k\ldots d_0)_{-3/2}:=\tfrac{1}{2}\sum_{i=0}^{k}d_i\left(-\tfrac32\right)^i \qquad (k\in\mathbb{N},\; d_i\in\mathcal{D})
\]	

\noindent and call $(d_k\ldots d_0)_{-3/2}$ an {\it integer ($-\frac32$)-expansion}. 
We proved in Proposition~\ref{prop:N} that each $N\in\Z$ has a unique integer ($-\frac32$)-expansion.
For instance, $-3=(2110)_{-3/2}$ and $4=(21122)_{-3/2}$. 

As a next step we characterize the set 
\[
\tfrac12\mathcal{D}[-\tfrac32]=\big\{(d_k\ldots d_0)_{-3/2}\;|\; k\in\mathbb{N},\; d_i\in\mathcal{D}\big\}
\]	
of all real numbers with an integer $(-\frac32)$-expansion (see also \cite[Example~3.3]{SSTW}).

\begin{theorem}\label{teo2}
The set of all numbers having an integer $(-\frac32)$-expansion is $\mathbb{Z}[\frac12]$. Here, as usual, we set $\mathbb{Z}[\frac12]=\{a2^{-\ell}\,\mid\, a\in\mathbb{Z},\,\ell\in\mathbb{N}\}$.
\end{theorem}

\begin{proof}
We have to show that $\tfrac12 \mathcal{D}[-\tfrac32]=\mathbb{Z}[\tfrac12].$
The inclusion $\tfrac12 \mathcal{D}[-\tfrac32] \subset \mathbb{Z}[\tfrac12]$ is trivial. For the reverse inclusion, let $N_0\in\mathbb{Z}[\tfrac12]$ be arbitrary. There exist $a\in\mathbb{Z}$ and $\ell\in\mathbb{N}$ such that $N_0=a2^{-\ell}$. We need to show that $N_0$ has an expansion of the form \eqref{eq:Nrep} for some $k\in \mathbb{N}$. We define the same recurrence as in \eqref{eq:algo} but forcing $N_i\in\mathbb{Z}[\tfrac12]$ for $i\geqslant0$, and proceed to show that $N_k=0$ for some $k\in\mathbb{N}$. 

Since $N_1$ is given by $2N_0=-3N_1+d_0$, where $N_1\in\mathbb{Z}[\tfrac12]$ and $d_0\in\mathcal{D}$, we have		
\begin{equation}\label{n1}
-3N_1=2N_0-d_0=\frac{a}{2^{\ell-1}}-d_0=\frac{a-2^{\ell-1}d_0}{2^{\ell-1}}.
\end{equation}	
To guarantee that $N_1\in\mathbb{Z}[\tfrac12]$, the numerator of this fraction has to be divisible by $3$, namely, we need to choose $d_0$ in a way that  $2^{\ell-1}d_0 \equiv a\pmod{3}$. As the inverse of $2^{\ell-1}$ in $\mathbb{Z}/3\mathbb{Z}$ is $2^{\ell-1}$, we get $d_0\equiv 2^{\ell-1}a\pmod{3}$ and $N_1$ is uniquely defined by~(\ref{n1}).
		
		
Iterating \eqref{eq:algo} yields that $N_i=a_i{2^{-\ell+i}}$ for some $a_i\in\mathbb{Z}$ for $i\in \{0,\ldots,\ell\}$. After $\ell$ steps, we get $N_\ell\in\mathbb{Z}$, and we are in the case covered by Proposition~\ref{prop:N}. This implies that there is $k_0\in \mathbb{N}$ such that $N_k=0$ for $k\geqslant k_0$ and thus~$N_0 \in \tfrac12 \mathcal{D}[-\tfrac32]$.
\end{proof}

By residue class considerations one can show that each $z\in\mathbb{Z}[\frac12]$ has a {\em unique} integer $(-\frac32)$-expansion. For instance, $-\frac38=(120)_{-3/2}$ and $\frac78=(111)_{-3/2}$.

\section{Expansions with fractional part: The reals.} \label{section 3}
We proceed to study $(-\frac32)$-expansions for arbitrary reals, allowing negative powers of the base. We will consider the matter of uniqueness and motivate the future construction of a tiling.

A desirable property of a number system is that almost all numbers (in a measure theoretic sense) can be expanded in a unique way. 
For example, although the decimal expansion is not always unique ({\it e.g.}\ the number $4$ can also be written as $3.9999...$), the set of numbers admitting more than one decimal expansion is very small in the sense that it has Lebesgue measure zero. This fact is reflected by the following {\em tiling property}. Consider the set of \textit{fractional parts} in the decimal number system, that is, the set of numbers that can be expanded using only negative powers of $10$: it corresponds to the unit interval $[0,1]$. On the other hand, the set of numbers whose expansion uses only nonnegative powers of $10$, that is, numbers with an \textit{integer expansion}, is equal to $\mathbb{Z}$ (if we permit the use of the minus sign). Consider the collection \begin{equation}\label{eq:C}
\{[0,1]+z\;|\;z\in\mathbb{Z}\}.
\end{equation}
This collection covers $\R$, and the only overlaps occur in the boundary points of the intervals, which form a measure zero set. We thus say that  $\{[0,1]+z\;|\;z\in\mathbb{Z}\}$ forms a \textit{tiling} of the real line. Here $[0,1]$ is the \textit{central tile} and $\mathbb{Z}$ is the \textit{translation set}. This tiling property is a geometric interpretation of the fact that almost all real numbers admit a unique expansion in the decimal system, and it works the same for any other $q$-ary number system ($q\in\mathbb{Z}$; $|q|\geqslant 2$).	
However, we will see that for a real number the representation in the number system $(-\frac32,\D)$ is {\it a priori} never unique. Let
\begin{equation}\label{expansion}
(d_k\ldots d_0.d_{-1}d_{-2}\ldots)_{-3/2}:=\tfrac{1}{2}\sum_{i=-\infty}^{k}d_i\left(-\tfrac32\right)^i \qquad (k\in\mathbb{Z},\; d_i\in\mathcal{D})
\end{equation}
and consider  
\begin{equation}\label{eq:Omegadef}
\varOmega=\{(0.d_{-1}d_{-2}\ldots)_{-3/2}\;|\; d_{-i}\in\mathcal{D}\},
\end{equation}
the set of \textit{fractional parts} in $(-\frac32,\D)$.
One can prove that $\varOmega= \big[-\tfrac{6}{5},\tfrac{4}{5}\big]$ by using the fact that $-\tfrac{3}2\varOmega=\varOmega\cup\big(\varOmega+\tfrac12\big)\cup(\varOmega+1)$ (we revisit this idea later on in Section~\ref{section 5}).  

Each decomposition of a real number $x$ as the sum of an element of $\mathbb{Z}[\tfrac12]$ and an element of $\varOmega$ leads to  an expansion of $x$ in the form (\ref{expansion}), by Theorem \ref{teo2} and the definition of $\varOmega$. Because the collection $\{\varOmega+z\,\mid\,z\in\mathbb{Z}[\tfrac12]\}$ covers the real line, each real number can be written as such a sum, and hence admits an expansion of the form  (\ref{expansion}). But since each $x\in\R$ is contained in multiple elements of the collection  $\{\varOmega+z\,\mid\,z\in\mathbb{Z}[\tfrac12]\}$ (in fact, in infinitely many), it admits multiple expansions of the form  (\ref{expansion}). For example, $\tfrac45=(0.020202\ldots)_{-3/2}=(2.11111\ldots)_{-3/2}$.  


Different translations of $\varOmega$ by elements of $\mathbb{Z}[\tfrac12]$ overlap in sets of positive measure; in other words, we do not have the desired tiling property. This results in expansions that are not unique. In the subsequent sections, we will find a way to ``embed'' the collection $\{\varOmega+z\,\mid\,z\in\mathbb{Z}[\tfrac12]\}$ in a suitable space where it will give rise to a tiling. 
	
	\section{The representation space.} \label{section 4}
	The real line seems to be ``too small" for the collection $\{\big[-\tfrac{6}{5},\tfrac{4}{5}\big]+z\,\mid\,z\in\mathbb{Z}[\tfrac12]\}$, so we wish to enlarge the space $\R$ in order to mend the issue with the overlaps\footnote{Another way, which we do not pursue here, would be to restrict the ``admissible'' digit strings, see \cite{MR2448050}.}. 
	Indeed, our next goal is to define a new space, called $\K$, in which the number system $(-\frac32,\D)$ induces a tiling in a natural way. The idea behind this is as follows: the overlaps occur because the three digits $\{0,1,2\}$ are ``too many'' for a base whose modulus is $\frac32$. Such a base would need one and a half digits, which is of course not doable. What causes all the problems is the denominator $2$. Roughly speaking, this denominator piles up powers of two which are responsible for the overlaps. It turns out that these overlaps can be  ``unfolded'' by adding a {\em $2$-adic factor} to  our representation space. 
The strategy of enlarging the representation space by $p$-adic factors that we are about to present was used in the setting of substitution dynamical systems {\it e.g.}~by Siegel~\cite{MR1997975} and in a much more general framework  than ours in~\cite{MR3391902}.
	
We begin by introducing the $2$-adic numbers; for more on this topic we refer the reader to \cite{MR1865659}.	Consider a nonzero rational number $y$ and write $y=2^\ell\tfrac{p}{q}$ 
where $\ell\in\mathbb{Z}$ and both $p$ and $q$ are odd. The $2$-adic absolute value in $\mathbb{Q}$ is defined by	
\[
\vert y \vert_2=  \begin{cases}
		2^{-\ell}, &   \hbox{if }  y \neq 0, \\
		 0, &  \hbox{if }  y=0,
	\end{cases}
\]	
and the $2$-adic distance between two rationals $x$ and $y$ is given by $|x-y|_2$.
Two points are close under this metric if their difference is divisible by a large positive power of~$2$. 

We define $\mathbb{Q}_2$ to be the completion of $\mathbb{Q}$ with respect to $|\cdot|_2$. The space $\mathbb{Q}_2$ is a field called the \textit{field of $2$-adic numbers}. Every nonzero $y\in\mathbb{Q}_2$ can be written uniquely as a series 
\[
y=\sum_{i=\ell}^{\infty}c_i2^i\qquad (\ell\in\mathbb{Z},\; c_i\in\{0,1\},\;c_\ell\neq0 ).
\]
This series converges in $\mathbb{Q}_2$ 
because large powers of two have small $2$-adic absolute value. Indeed, we have $|y|_2=2^{-\ell}$.
	
We define our \textit{representation space} as $\K=\mathbb{R}\times\mathbb{Q}_2$, with the  additive group structure given by componentwise addition. Moreover, $\mathbb{Z}[\tfrac12]$ acts on $\K$ by multiplication, more precisely, if $\alpha\in\mathbb{Z}[\tfrac12]$ and $(x_1,x_2)\in\K$ then
\[
\alpha\cdot(x_1,x_2)=(\alpha x_1,\alpha x_2)=(x_1,x_2)\cdot\alpha.
\]

 For every $(x_1,x_2),(y_1,y_2)\in\K$ define 
\[
\textbf{d}((x_1,x_2),(y_1,y_2)):=\max\{|x_1-y_1|,|x_2-y_2|_2\}.
\]
Then $\textbf{d}$ is a metric on $\K$. Intuitively, two points in $\K$ are far apart if either their real components are far apart or their $2$-adic components are far apart.
	
We define the embedding	
	\begin{equation*}
	\begin{split}
	\varphi: \,\mathbb{Q}\to \K,\quad z\mapsto(z,z).
	\end{split}
	\end{equation*}		
Consider the image of $\mathbb{Z}[\tfrac12]$ under $\varphi$. Despite both coordinates of $\varphi(z)$ being the same, it does not lie in a diagonal. Indeed, points of $\mathbb{Z}[\tfrac12]$ that are close in $\R$ are far apart in the $2$-adic distance. In particular, we will show that the points of $\varphi(\mathbb{Z}[\tfrac12])$ form a \emph{lattice}.

A subset $\Lambda$ of $\K$ is a {\em lattice} if it satisfies the three following conditions.
	\begin{enumerate}
	 	\item $\Lambda$ is a group.
	 	 \item $\Lambda$ is uniformly discrete, meaning there exists $r>0$ such that every open ball of radius $r$ in $\K$ contains at most one point of $\Lambda$.
	 	 \item $\Lambda$ is relatively dense, meaning there exists $R>0$ such that every closed ball of radius $R$ in $\K$ contains at least one point of $\Lambda$.
	\end{enumerate}
	 	 
	\begin{lemma}\label{lem:lattice}
		 $\varphi(\mathbb{Z}[\tfrac12])$ is a lattice in $\K$.
	\end{lemma}
	 
	 \begin{proof}
	 	
	 \begin{enumerate}
	 	\item  The fact that $\varphi(\mathbb{Z}[\tfrac12])$ is a group follows from the additive group structure of $\mathbb{Z}[\tfrac12]$ because $\varphi$ is a group homomorphism.
	 
	 \item To get uniform discreteness of $\varphi(\mathbb{Z}[\tfrac12])$ we show first that $\textbf{d}(\varphi(z),\varphi(0))\geqslant1$ holds for every nonzero $z\in\mathbb{Z}[\tfrac12]$. Recall that $\textbf{d}(\varphi(z),\varphi(0))=\max\{|z|,|z|_2\}$. If $|z|<1$ there exist $a,\ell\in\mathbb{Z}$ with $a$ odd and $\ell\ge1$ such that $z=a2^{-\ell}$, so $|z|_2=2^\ell> 1$ and thus $\textbf{d}(\varphi(z),\varphi(0))\geqslant1$. Because of the group structure, this implies that the distance between any two elements of $\varphi(\mathbb{Z}[\tfrac12])$ is at least one, hence $\varphi(\mathbb{Z}[\tfrac12])$ is uniformly discrete.
	 
	 \item For the relative denseness, consider an arbitrary element $(x_1,x_2)\in\K$. We claim that there exists $z\in\mathbb{Z}[\tfrac12]$ such that $\textbf{d}(\varphi(z),(x_1,x_2))\leqslant2$.
	 Let $z_1\in\mathbb{Z}$ be one of the integers being closest to $x_1$. If $x_2=\sum_{i=\ell}^{\infty}c_i2^i$ with $\ell\in\mathbb{Z}$ and $c_i\in\{0,1\}$ then set $z_2=\sum_{i=\ell}^{-1}c_i2^i\in\mathbb{Z}[\tfrac12]$ (note that $z_2=0$ if $\ell \geqslant 0$). Therefore, $$\textbf{d}((x_1,x_2),(z_1,z_2))=\max\{|x_1-z_1|,|x_2-z_2|_2\}\leqslant1.$$ Now we set $z=z_1+z_2\in\mathbb{Z}[\tfrac12]$. Because $z_1$ is an integer, $|z_1|_2\leqslant 1$, and since $z_2\in[0,1]$, we have $|z_2|\leqslant 1$.
	 Thus $\textbf{d}(\varphi(z),(z_1,z_2))=\max\{|z_2|,|z_1|_2\}\leqslant1$, and so $\textbf{d}(\varphi(z),(x_1,x_2))\leqslant2$ by the triangle inequality, hence $\varphi(\mathbb{Z}[\tfrac12])$ is relatively dense. \qedhere
	\end{enumerate}
	\end{proof}
	
\begin{figure*}[h!]
		\centering\includegraphics[scale=0.6]{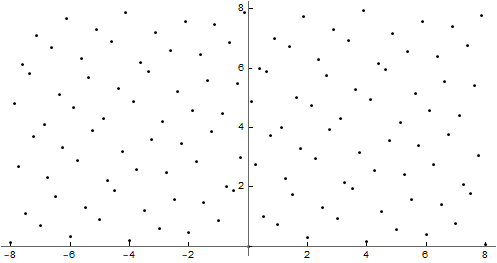}
		\caption{Representation of the lattice points $\varphi(\frac j{8})$ for $-64\leqslant j\leqslant 64$ in $\R^2$.}
		\label{points}
	\end{figure*}

Figure~\ref{points} illustrates some points of $\varphi(\mathbb{Z}[\tfrac12])$. Drawing pictures in this setting is not straightforward: the space $\K$ is non-Euclidean, so we need to represent it in~$\R^2$ while somehow maintaining the $2$-adic nature of the second component. We do this in the following way: any point $y\in\Q_2$ can be written uniquely (up to ``leading zeros'') as a series $y=\frac12\sum_{i=\ell}^{\infty}c_i(-\tfrac23)^{i}$ with $\ell\in\mathbb{Z},\; c_i\in\{0,1\}$, that converges in the $2$-adic metric (this is a $2$-adic expansion, not a $(-\tfrac32)$-expansion!). We consider the mapping
\begin{equation}\label{gamma}
\gamma:\mathbb{Q}_2\to \mathbb{R};\qquad \tfrac12\sum_{i=\ell}^{\infty}c_i(-\tfrac23)^{i}\mapsto\sum_{i=\ell}^{\infty}c_i2^{-i},
\end{equation} 
which is well defined since the sum on the right hand side converges in $\R$. A point $(x_1,x_2)\in \K$ is now represented as $(x_1,\gamma(x_2))\in \R^2$.

The lattice $\varphi(\mathbb{Z}[\tfrac12])$, which will play the role of the ``integers'' in $\mathbb{K}$, turns out to be a proper translation set for a tiling of $\K$ related to the number system $(-\tfrac32,\mathcal{D})$. In the real case, when defining a tiling we allowed overlaps as long as it was on a set of Lebesgue measure zero. In order to generalize this, we need to define a natural measure on $\K$.

A \textit{Haar measure} is a translation invariant Borel measure that is finite for compact sets. It can be defined in spaces with a sufficiently ``nice" structure (more specifically, it is defined on locally compact topological groups). Such a measure is unique up to a scaling factor. The Lebesgue measure $\mu_\infty$ in $\mathbb{R}$ is a Haar measure. 

Let $\mu_2$ be the Haar measure in $\mathbb{Q}_2$ that satisfies 
\[
\mu_2(2^\ell\mathbb{Z}_2)=2^{-\ell},
\]	
where $\mathbb{Z}_2=\big\{\sum_{i=0}^{\infty}c_i2^i\;|\;c_i\in\{0,1\} \big\}\subset \Q_2$ is the ring of $2$-adic integers. This is a very natural measure: multiplying by large powers of two makes a set small in measure. 

Let $\mu=\mu_\infty\times\mu_2$ be the product measure of $\mu_\infty$ and $\mu_2$ on $\K=\mathbb{R} \times \mathbb{Q}_2$, that is, if $M_1\subset\mathbb{R}$ and $M_2\subset\mathbb{Q}_2$ are respectively measurable, then the sets of the form $M=M_1\times M_2$ generate the $\sigma$-algebra of $\mu$, and $\mu(M):=\mu_\infty(M_1)\mu_2(M_2)$. One can show that $\mu$ is a Haar measure on~$\K$.

	For any measurable set $M=M_1\times M_2\subset\K$, we have
\begin{equation*}\label{muinfty}
\mu_\infty(-\tfrac32M_1)=\tfrac{3}{2}\mu_\infty(M_1),\qquad \mu_2(-\tfrac32M_2)=2\mu_2(M_2)
\end{equation*}		

\noindent	which yields
	\begin{equation}\label{eq:muexpand}
	\mu(-\tfrac32M)=\mu_\infty(-\tfrac32M_1)\mu_2(-\tfrac32M_2)=3\mu(M).
	\end{equation}

\noindent Thus multiplying any measurable set $M\subset\K$ by the base $-\frac32$ enlarges the measure by the factor $3$, which can be interpreted as having ``enough space'' for three digits.
	
\section{The tile $\mathcal{F}$.} \label{section 5}
In this section we define a set $\mathcal{F}\subset\K$ that plays the same role for $(-\frac32,\D)$ as the unit interval does for the decimal system. We explore some of its topological and measure theoretic properties.

Recall that in \eqref{eq:Omegadef} we defined the set $\varOmega$ of fractional parts, consisting of elements of the form $(0.d_{-1}d_{-2}\ldots)_{-3/2}$. We now embed the digits in $\K$, obtaining the set
\[
\mathcal{F}:=\Big\{\tfrac12\sum_{i=1}^{\infty}\varphi(d_{-i})\left(-\tfrac32\right)^{-i}\;|\;d_{-i}\in\D\Big\}.
\]
The set $\mathcal{F}$ is a compact subset of $\K$. Indeed, given any sequence in $\mathcal{F}$, we use a Cantor diagonal argument to find a convergent subsequence.

 Let $x\in\mathcal{F}$: if we multiply $x$ by the base $-\tfrac{3}{2}$, we obtain $-\tfrac{3}{2}x\in\mathcal{F}+\tfrac12\varphi(d_{-1})$ with $d_{-1}\in\{0,1,2\}$ (this can be interpreted as the analog of moving the decimal point one place to the right). Thus $\mathcal{F}$ satisfies the \textit{set equation}
\begin{equation}\label{eq:seteqF}
-\tfrac{3}{2}\mathcal{F}=\mathcal{F}\cup\big(\mathcal{F}+\varphi(\tfrac12)\big)\cup\big(\mathcal{F}+\varphi(1)\big)
\end{equation}
in $\K$, which can be written shortly as
$
-\tfrac32\mathcal{F}=\mathcal{F}+\tfrac12\varphi(\mathcal{D}).
$
It turns out that this set equation completely characterizes $\mathcal{F}$. Note that (\ref{eq:seteqF}) is equivalent to 
\begin{equation}\label{eq:seteqF2}
\mathcal{F}=\left( -\tfrac23 \right)\mathcal{F}\cup\left( -\tfrac23 \right)\big(\mathcal{F}+\varphi(\tfrac12)\big)\cup\left( -\tfrac23 \right)\big(\mathcal{F}+\varphi(1)\big),
\end{equation}
and multiplying by $-\tfrac23$ is a uniform contraction in $\K$: it is a contraction in $\mathbb{R}$ because $\vert-\tfrac23\vert<1$ and also in $\mathbb{Q}_2$ because $|-\tfrac23|_2=\tfrac12<1.$
Thus \eqref{eq:seteqF2} states that $\mathcal{F}$ is equal to the union of three contracted copies of itself. Because of this contraction property we may apply Hutchinson's Theorem (see \cite{MR625600}) which says that there exists a unique nonempty compact subset of $\K$ that satisfies the set equation \eqref{eq:seteqF2}. Thus $\mathcal{F}$ is uniquely defined as the nonempty compact set satisfying (\ref{eq:seteqF2}) (or, equivalently, (\ref{eq:seteqF})). The set $\mathcal{F}$ is called a {\em rational self-affine tile} in the sense of~\cite{MR3391902}.  
	
Since according to Theorem~\ref{teo2} the set $\mathbb{Z}[\frac12]$ is the analog of $\Z$ in the number system $(-\frac32,\D)$, we define the analog of the collection in \eqref{eq:C} by setting
\[
\mathcal{C}=\{\mathcal{F}+\varphi(z)\,\mid\, z\in \mathbb{Z}[\tfrac12] \}.
\]
Then $\mathcal{C}$ is a collection of copies of $\mathcal{F}$ translated by elements of the lattice $\varphi(\mathbb{Z}[\tfrac12])$. We will show in Theorem~\ref{th:tiling} that $\mathcal{C}$ is a \emph{tiling} of $\K$, meaning that:
	\begin{enumerate}
		\item $\mathcal{C}$ is a {\em covering} of $\K$, {\it i.e.}, $\langle \mathcal{C}\rangle=\K$, where $\langle \mathcal{C}\rangle=\mathcal{F}+\varphi(\mathbb{Z}[\tfrac12])
		$ is the union of the elements of $\mathcal{C}$.
		\item Almost every point in $\K$ (with respect to the measure $\mu$) is contained in exactly one element of $\mathcal{C}$.
	\end{enumerate}

Figure~\ref{tile} shows a representation of $\mathcal{F}$ in~$\R^2$, again using the function $\gamma$ from (\ref{gamma}) to map $\mathbb{Q}_2$ to $\mathbb{R}$. 
\begin{figure*}[h!]
\centering\includegraphics*[width=0.7\linewidth]{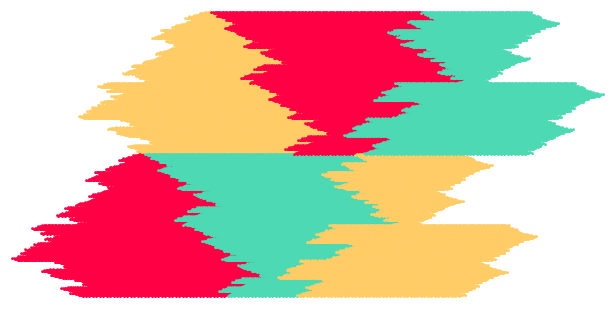}
\caption{A patch of the tiling $\mathcal{C}$ of $\K$ by translates of $\mathcal{F}$.}
\label{patch}
\end{figure*}
Figure~\ref{patch} shows a patch of $\mathcal{C}$; the translates of $\mathcal{F}$ appear to have different shapes, but this is due to the embedding of $\K$ in $\R^2$. The illustration indicates, however, that the different translates of $\mathcal{F}$ do not overlap other than in their boundaries.

Because $\mathcal{F}$ plays the role of the central tile, we require it to be a reasonably nice set topologically speaking. As a next step, we prove that $\mathcal{F}$ is the closure of its interior and that its boundary $\partial \mathcal{F}$ has measure zero. 
In a general setting, this result is contained in \cite[Theorem~1]{MR3391902}.

\begin{theorem} \label{th:int}
$\mathcal{F}$ is the closure of its interior.
\end{theorem}

\begin{proof}

We first prove that $\mathcal{C}$ is a covering of $\K$, {\it i.e.}, $\langle \mathcal{C}\rangle=\K$, where $\langle \mathcal{C}\rangle$ is the union of the elements of $\mathcal{C}$. Applying the set equation (\ref{eq:seteqF})
we obtain 	
\begin{equation*}
-\tfrac32\langle \mathcal{C}\rangle=-\tfrac32\mathcal{F}-\tfrac32\varphi(\mathbb{Z}[\tfrac12])=\mathcal{F}+\tfrac12\varphi(\mathcal{D})-\tfrac32\varphi(\mathbb{Z}[\tfrac12]).
\end{equation*}	
Note that $\tfrac12\mathcal{D}=\{0,\tfrac12,1\}$ is a complete set of representatives of residue classes of $\mathbb{Z}[\tfrac12]/(-\tfrac32)\mathbb{Z}[\tfrac12]$, so 
\[
\tfrac12\varphi(\mathcal{D})-\tfrac32\varphi(\mathbb{Z}[\tfrac12])=\varphi(\mathbb{Z}[\tfrac12]).
\]	
Thus $-\tfrac32 \langle \mathcal{C}\rangle=\langle \mathcal{C}\rangle$ and, {\it a fortiori}, for any $k\in\mathbb{N}$ we have
$(-\tfrac23)^k\langle \mathcal{C}\rangle=\langle \mathcal{C}\rangle$.
Recall that multiplying by $-\tfrac23$ is a contraction in $\K$. 		
We have shown in Lemma~\ref{lem:lattice} that $\varphi(\mathbb{Z}[\tfrac12])$ is a relatively dense set in $\K$, and therefore so is $\langle \mathcal{C}\rangle$, meaning there is some $R>0$ for which every closed ball of radius $R$ intersects $\langle \mathcal{C}\rangle$. But since $\langle \mathcal{C}\rangle$ is invariant under contractions by $(-\tfrac23)^k$, this implies that any ball of radius  $(\frac23)^kR$ with $k\in\mathbb{N}$ intersects $\langle \mathcal{C}\rangle$, hence it is dense in $\K$.

Consider now an arbitrary point in $x\in\K$ and a bounded neighborhood $V$ of $x$. Since $\varphi(\mathbb{Z}[\tfrac12])$ is uniformly discrete and $V$ is bounded, $V$ intersects only a finite number of translates of $\mathcal{F}$, each of which is compact. Since $\langle \mathcal{C}\rangle$ is dense in $\K$, $x$ cannot be at positive distance from all these translates of $\mathcal{F}$. Thus $x$ is contained in some translate of $\mathcal{F}$ and, hence, $x\in\langle \mathcal{C}\rangle$. Since $x$ was arbitrary this implies that $\langle \mathcal{C}\rangle=\K$.
		
Next, we show that $\mbox{int}\,\mathcal{F}\neq\varnothing$. Assume on the contrary that $\mbox{int}\,\mathcal{F}=\varnothing$. Consider the sets
\[
U_z:=\K\setminus(\mathcal{F}+\varphi(z)) \qquad(z\in\mathbb{Z}[\tfrac12]).
\]
By assumption, $U_z$ is dense in $\K$ for each $z\in\mathbb{Z}[\tfrac12]$, and $\{U_z\;\mid \, z\in\mathbb{Z}[\tfrac12]\}$ is a countable collection. Baire's theorem asserts that a countable intersection of dense sets is dense. 
But 
\[
\bigcap_{z\in\mathbb{Z}[\frac12]}U_z=\K\setminus\bigcup_{z\in\mathbb{Z}[\frac12]}\mathcal{F}+\varphi(z)=\K\setminus\langle\mathcal{C}\rangle=\varnothing,
\]
which is clearly not dense. This contradiction yields $\mbox{int}\,\mathcal{F}\neq\varnothing$.

We now prove the result. Iterating the set equation \eqref{eq:seteqF} for $k\in\mathbb{N}$ times yields 
\[
\mathcal{F}=\left(-\tfrac23\right)^k\mathcal{F}+\tfrac12\left(\varphi(\mathcal{D})\left(-\tfrac23\right)^k+\varphi(\mathcal{D})\left(-\tfrac23\right)^{k-1}+\cdots+\varphi(\mathcal{D})\left(-\tfrac23\right)\right).
\]		
Setting 
\begin{equation}\label{eq:Dk}
\mathcal{D}_k:=\mathcal{D}+\mathcal{D}\left(-\tfrac32\right)+\cdots+\mathcal{D}\left(-\tfrac32\right)^{k-1}\end{equation}
this becomes
\begin{equation}\label{eq:itset}
\mathcal{F}=(-\tfrac23)^k(\mathcal{F}+\tfrac12\varphi(\mathcal{D}_k))\qquad	
(k\in\mathbb{N}),
\end{equation}  which means we can write $\mathcal{F}$ as a finite union of arbitrarily small shrunk translated copies of itself. We know that $\mathcal{F}$ has an inner point $x$, therefore each copy of the form $(-\tfrac23)^k(\mathcal{F}+\tfrac12\varphi(d)),\,d\in\D_k$, has an inner point. Thus for any $y\in\mathcal{F}$ and any $\varepsilon>0$ we can choose $k\in\mathbb{N}$ and $d\in \mathcal{D}_k$ so that $\mathrm{diam} ((-\tfrac23)^k(\mathcal{F}+\tfrac12\varphi(d))) <  \varepsilon$ and $y \in (-\tfrac23)^k(\mathcal{F}+\tfrac12\varphi(d))$. Thus there is an inner point at distance less than $\varepsilon$ from $y$. Since  $y\in\mathcal{F}$ and $\varepsilon>0$ were arbitrary, this proves that $\mathcal{F}$ is the closure of its interior.
\end{proof}

\begin{theorem}\label{th:zeromeasure}
The boundary of $\mathcal{F}$ has measure zero.
\end{theorem}

\begin{proof}
Let $x$ be an inner point of $\mathcal{F}$ and $B_\varepsilon(x)\subset\mathcal{F}$ an open ball of radius $\varepsilon>0$ centered at $x$. Because multiplication by $-\tfrac23$ is a uniform contraction in $\mathbb{K}$, there is $k\in\mathbb{N}$ such that
$\mathrm{diam}\,(-\tfrac23)^k\mathcal{F}<\varepsilon$. Thus by (\ref{eq:itset}) there is $d_0\in\mathcal{D}_k$ such that	
\[
\left(-\tfrac23\right)^k\big(\mathcal{F}+\tfrac12\varphi(d_0)\big)\subset B_\varepsilon(x)\subset\mbox{int}\,\mathcal{F}.
\]	
Let 
\[
y\in\partial\big(\left(-\tfrac23\right)^k\big(\mathcal{F}+\tfrac12\varphi(d_0)\big)\big)\subset\mbox{int}\,\mathcal{F}.
\]	
Since $y$ is also an inner point of $\mathcal{F}$, and (\ref{eq:itset}) writes $\mathcal{F}$ as a finite union of compact sets, $y$ must necessarily lie in $(-\tfrac23)^k(\mathcal{F}+\tfrac12\varphi(d))$ for some $d\in \mathcal{D}_k\setminus\{d_0\}$.  Thus the boundary $\partial((-\tfrac23)^k(\mathcal{F}+\tfrac12\varphi(d_0)))$ is covered at least twice by the collection $\{(-\tfrac23)^k(\mathcal{F}+\tfrac12\varphi(d))\;|\;d\in\mathcal{D}_k\}$. This entails that
\begin{equation*}
		\begin{split}
		\mu(\mathcal{F})&=\mu\big(\bigcup_{d\in\mathcal{D}_k}\left(-\tfrac23\right)^k\big(\mathcal{F}+\tfrac12\varphi(d)\big)\big)\\
		&\leqslant\sum_{d\in\mathcal{D}_k}\mu\big(\left(-\tfrac23\right)^k\big(\mathcal{F}+\tfrac12\varphi(d)\big)\big)-\mu\big(\partial\big(\left(-\tfrac23\right)^k\big(\mathcal{F}+\tfrac12\varphi(d_0)\big)\big)\big).
		\end{split}
\end{equation*}		
Note that (as a Haar measure) $\mu$ is translation invariant, the cardinality of $\D_k$ is $3^{k}$, and from (\ref{eq:muexpand}) it follows that $\mu((-\tfrac23)^k\mathcal{F})=3^{-k}\mu(\mathcal{F})$. All this combined yields		
\[
\begin{split}
\mu(\mathcal{F})&\leqslant\sum_{d\in\mathcal{D}_k}\mu\big(\left(-\tfrac23\right)^k\mathcal{F}\big)-\mu\big(\partial\big(\left(-\tfrac23\right)^k\mathcal{F}\big)\big)
\\&\leqslant3^k3^{-k}\mu(\mathcal{F})-\mu\big(\partial\big(\left(-\tfrac23\right)^k\mathcal{F}\big)\big)\\
&=\mu(\mathcal{F})-\mu\big(\partial\big(\left(-\tfrac23\right)^k\mathcal{F}\big)\big)	
\end{split}		
\]
and therefore $\mu(\partial((-\tfrac23)^k\mathcal{F}))=0$. This implies that $\mu(\partial\mathcal{F})=0$.	
\end{proof}


\section{The tiling.}	\label{section 7}

This section contains our final result: a tiling theorem for the $(-\frac32)$-number system. This result is contained in \cite[Theorem~2]{MR3391902} in a more general setting. For our special case, the proof is much simpler. 
As mentioned before, the tiling property is important because it relates to the uniqueness (almost everywhere) of expansions in the $(-\frac32)$-number system embedded in $\K$. We prove this as a corollary of our tiling theorem.

\begin{theorem}\label{th:tiling}
The collection $\mathcal{C}=\{\mathcal{F}+\varphi(z)\,\mid\, z\in \mathbb{Z}[\tfrac12] \}$ forms a tiling of~$\,\K$.	 
\end{theorem}

\begin{proof}
We have shown in the proof of Theorem \ref{th:int} that $\mathcal{C}$ is a covering of $\K$. It remains to show that almost every point of $\K$ is covered by exactly one element of the collection~$\mathcal{C}$.
Recall that for each $k\geqslant1$, the sets $\tfrac12\D_k$ (see (\ref{eq:Dk})) consist of all the integer ($-\tfrac32$)-expansions with at most $k$ digits. According to Theorem~\ref{teo2}, the set $\Z[\tfrac12]$ is the set of all integer ($-\tfrac32$)-expansions. This implies that $\Z[\tfrac12]=\bigcup_{k\geqslant 1}\tfrac12\D_k$ and, hence,
$
\K=\mathcal{F}+\varphi(\Z[\tfrac12])=\bigcup_{k\geqslant 1}\mathcal{F}+\tfrac12\varphi(\D_k).
$
Therefore, it suffices to prove that the collection $\{\mathcal{F}+\tfrac12\varphi(d)\,\mid\,d\in\mathcal{D}_k\}$ has \textit{essentially disjoint} elements for each $k\geqslant1$, that is, if $d,d'\in\D_k$ are distinct then 
$\mu((\mathcal{F}+\tfrac12\varphi(d))\cap(\mathcal{F}+\tfrac12\varphi(d')))=0$.
Applying (\ref{eq:itset}) we obtain
\begin{equation*}
\begin{split}	
3^k\mu(\mathcal{F})&=\mu((-\tfrac32)^k\mathcal{F})=\mu(\bigcup_{d\in\D_k}\mathcal{F}+\tfrac12\varphi(d))\leqslant \sum_{d\in\D_k}\mu(\mathcal{F}+\tfrac12\varphi(d))=3^k\mu(\mathcal{F}).
\end{split}
\end{equation*}
This implies equality everywhere and, hence, different $\tfrac12\varphi(\D_k)$-translates of $\mathcal{F}$ only overlap in sets of measure zero. Thus the same is true for different $\Z[\tfrac12]$-translates of $\mathcal{F}$. So the tiles in $\mathcal{C}$ are essentially disjoint, and $\mathcal{C}$ is a tiling.	
\end{proof}

\begin{corollary}\label{cor:unique}
Almost every point $x\in \K$ has a unique expansion of the form
\begin{equation}\label{eq:phirep}
x=\tfrac12\sum_{i=-\infty}^k \left(-\tfrac32\right)^i\varphi(d_i) \qquad (k\in\mathbb{N},\; d_i\in\D; \,d_k\neq0 \mbox{ whenever }k\geqslant 1).
\end{equation}
\end{corollary}

\begin{proof}	
Let $x\in\K$ and suppose it has two different expansions
\[x=\tfrac12\sum_{i=-\infty}^k \left(-\tfrac32\right)^{i}\varphi(d_{i})=\tfrac12\sum_{i=-\infty}^k \left(-\tfrac32\right)^{i}\varphi(d'_{i}),
\]
where $d_k\neq0$ for $k\ge 1$ and where we pad the second expansion with zeros if necessary. Let $m\leqslant k$ be the largest integer such that $d_{m}\neq d'_{m}$, and consider the point $(-\tfrac32)^{-m}x$. Recall that multiplying $x$ by $(-\tfrac32)^{-m}$ is the analog of moving the decimal point $m$ places to the left if $m$ is positive and to the right if it is negative. Let $\omega:=(d_{k}\ldots d_{m+1} d_{m})_{-3/2}$ and $\omega':=(d'_{k}\ldots d'_{m+1}d'_m)_{-3/2}$. Then $\omega,\omega'\in \Z[\tfrac12]$ are distinct, and it follows from our assumption and the definition of the tile $\mathcal{F}$ that $(-\tfrac32)^{-m}x-\varphi(\omega),(-\tfrac32)^{-m}x-\varphi(\omega')\in\mathcal{F}$. Hence, we obtain
$$(-\tfrac32)^{-m}x\in(\mathcal{F}+\varphi(\omega))\cap(\mathcal{F}+\varphi(\omega')).
$$
As tiles only overlap on their boundaries, this implies that $x\in(-\tfrac32)^{m}\partial(\mathcal{F}+\varphi(\omega))$. 
Therefore, a point $x\in\K$ has two different expansions if and only if $x\in\varGamma$, where $\varGamma:=\bigcup_{m\in\Z}(-\tfrac32)^{m}\partial(\mathcal{F}+\varphi(\Z[\tfrac12]))$.
Since $\Z[\tfrac12]$ is countable, $\varGamma$ is a countable union of the sets $(-\tfrac32)^{m}\partial(\mathcal{F}+\varphi(z))$, $m\in\mathbb{Z}$, $z\in\Z[\tfrac12]$, each of which has measure $0$. Thus $\mu(\varGamma)=0$, which gives the result.
\end{proof}

\bibliographystyle{plain}	
\bibliography{biblio}

\begin{thebibliography}{10}

\bibitem{MR2448050}
Shigeki Akiyama, Christiane Frougny, and Jacques Sakarovitch.
\newblock Powers of rationals modulo 1 and rational base number systems.
\newblock {\em Israel J. Math.}, 168:53--91, 2008.

\bibitem{MR3051451}
Petr Ambro\v{z}, Daniel Dombek, Zuzana Mas\'{a}kov\'{a}, and Edita
  Pelantov\'{a}.
\newblock Numbers with integer expansion in the numeration system with negative
  base.
\newblock {\em Funct. Approx. Comment. Math.}, 47(part 2):241--266, 2012.

\bibitem{MR2735753}
Val\'{e}rie Berth\'{e}, Anne Siegel, Wolfgang Steiner, Paul Surer, and
  J\"{o}rg~M. Thuswaldner.
\newblock Fractal tiles associated with shift radix systems.
\newblock {\em Adv. Math.}, 226(1):139--175, 2011.

\bibitem{MR1613189}
Daniel~M. Gordon.
\newblock A survey of fast exponentiation methods.
\newblock {\em J. Algorithms}, 27(1):129--146, 1998.

\bibitem{Gru}
Vittorio {Gr\"unwald}.
\newblock Intorno all'aritmetica dei sistemi numerici a base negativa con
  particolare riguardo al sistema numerico a base negativo-decimale per lo
  studio delle sue analogie coll'aritmetica (decimale).
\newblock {\em Giornale di Matematiche di Battaglini}, 23:203--221, 1885.
\newblock Errata, p. 367.

\bibitem{MR1865659}
Jan~E. Holly.
\newblock Pictures of ultrametric spaces, the {$p$}-adic numbers, and valued
  fields.
\newblock {\em Amer. Math. Monthly}, 108(8):721--728, 2001.

\bibitem{MR625600}
John~E. Hutchinson.
\newblock Fractals and self-similarity.
\newblock {\em Indiana Univ. Math. J.}, 30(5):713--747, 1981.

\bibitem{Kempner}
Aubrey~J. Kempner.
\newblock Anormal {S}ystems of {N}umeration.
\newblock {\em Amer. Math. Monthly}, 43(10):610--617, 1936.

\bibitem{Knuth:60}
Donald~E. Knuth.
\newblock An imaginary number system.
\newblock {\em Comm. ACM}, 3:245--247, 1960.

\bibitem{Knu}
Donald~E. Knuth.
\newblock {\em The art of computer programming. {V}ol. 2}.
\newblock Addison-Wesley, Reading, MA, 1998.
\newblock Seminumerical algorithms, Third edition.

\bibitem{MR2042759}
Susan Landau.
\newblock Polynomials in the nation's service: using algebra to design the
  advanced encryption standard.
\newblock {\em Amer. Math. Monthly}, 111(2):89--117, 2004.

\bibitem{Renyi:57}
Alfr\'{e}d R\'{e}nyi.
\newblock Representations for real numbers and their ergodic properties.
\newblock {\em Acta Math. Acad. Sci. Hungar.}, 8:477--493, 1957.

\bibitem{SSTW}
Klaus Scheicher, Paul Surer, J\"org~M. Thuswaldner, and Christiaan~E. van~de
  Woestijne.
\newblock Digit systems over commutative rings.
\newblock {\em Intern. J. Number Theory}, 10(6):1459--1483, 2014.

\bibitem{MR2024754}
Nikita Sidorov.
\newblock Almost every number has a continuum of {$\beta$}-expansions.
\newblock {\em Amer. Math. Monthly}, 110(9):838--842, 2003.

\bibitem{MR1997975}
Anne Siegel.
\newblock Repr\'{e}sentation des syst\`emes dynamiques substitutifs non
  unimodulaires.
\newblock {\em Ergodic Theory Dynam. Systems}, 23(4):1247--1273, 2003.

\bibitem{MR3391902}
Wolfgang Steiner and J\"{o}rg~M. Thuswaldner.
\newblock Rational self-affine tiles.
\newblock {\em Trans. Amer. Math. Soc.}, 367(11):7863--7894, 2015.

\end{thebibliography}

\end{document}